%% file: paper.tex
\documentclass[12pt]{article}

\usepackage[english]{babel}

\usepackage{amsfonts, amssymb, amsmath, epsf}
\usepackage{amsthm}

\usepackage{cite}
\usepackage{float}
\usepackage{graphicx}

\usepackage{fullpage}

\theoremstyle{plain}
\newtheorem{theorem}{Theorem}[section]
\newtheorem{lemma}{Lemma}[section]
\newtheorem{proposition}{Proposition}[section]

\newtheorem{condition}{Condition}[section]

\newtheorem{prp}{Proposition}[section]

\theoremstyle{definition}
\newtheorem{definition}{Definition}[section]

\newtheorem{remark}{Remark}[section]

\numberwithin{equation}{section} \numberwithin{figure}{section}

\renewcommand{\phi}{{\varphi}}

\newcommand{\cK}{{\mathcal K}}

\newcommand{\cP}{{\mathcal P}}

\newcommand{\cH}{{\mathcal H}}
\newcommand{\cS}{{\mathcal S}}

\newcommand{\ZZ}{\mathbb{Z}}

\newcommand{\bbR}{{\mathbb R}}
\newcommand{\bbQ}{{\mathbb Q}}
\newcommand{\bbN}{{\mathbb N}}

\newcommand{\bbZ}{{\mathbb Z}}

\newcommand{\Gr}{{\mathcal Y}}

\let\phi=\varphi

\newcommand{\ep}{\varepsilon}
\newcommand{\al}{\alpha}

\title{Rattling in spatially discrete diffusion equations with hysteresis}%[Spatially discrete diffusion equations with hysteresis]

\author{Pavel  Gurevich\thanks{Free University of Berlin, Institute of Mathematics I,
	Arnimallee 3, 14195, Berlin Germany; Peoples' Friendship University of Russia 117198, Moscow Miklukho-Maklaya str.~6, email: gurevich@math.fu-berlin.de}
\and
Sergey  Tikhomirov\thanks{Max Planck Institute for Mathematics in the Sciences,
	Inselstra\ss e 22, 04103 Leipzig, Germany; Chebyshev Laboratory, St. Petersburg State University, 14th Line, 29b, Saint Petersburg, 199178 Russia; email: sergey.tikhomirov@mis.mpg.de}}

\begin{document}

\maketitle
%\slugger{mms}{xxxx}{xx}{x}{x--x}

\begin{abstract}
  The paper treats a reaction-diffusion equation with hysteretic nonlinearity on a one-dimensional lattice. It arises as a result of the spatial discretization of the corresponding continuous model with so-called nontransverse initial data and exhibits a propagating microstructure --- which we call {\em rattling} --- in the hysteretic component of the solution. We analyze this microstructure and determine the speed of its propagation depending on the parameters of hysteresis and the nontransversality coefficient in the initial data.
\end{abstract}

\textbf{Key words.}
Spatially discrete parabolic equations, reaction-diffusion equations, lattice, hysteresis, pattern, rattling.
%\end{keywords}

\textbf{AMS subject classification.}
34K31, %	Lattice functional-differential equations
47J40, %Equations with hysteresis operators
35B36, %Pattern formation (PDE)
37L60.% Dynamical systems and ergodic theory/Infinite-dimensional dissipative dynamical systems/Lattice dynamics
%\end{AMS}

%\pagestyle{myheadings}
%\thispagestyle{plain}
%\markboth{PAVEL  GUREVICH, SERGEY  TIKHOMIROV}{RATTLING IN SPATIALLY DISCRETE %DIFFUSION EQUATIONS WITH HYSTERESIS}

%\tableofcontents

\input intro.tex

\input sect2.tex

\input sect3.tex

\input example.tex

\input sect4.tex

\input sect5.tex

\input sect6.tex

\input sect7.tex

\input appendixA.tex

\input bibl.tex
\end{document}

%% file: intro.tex
\section{Introduction}

\begin{figure}[ht]
\begin{center}
\begin{minipage}[b]{0.4\textwidth}
        \centering
        \includegraphics[width=4cm]{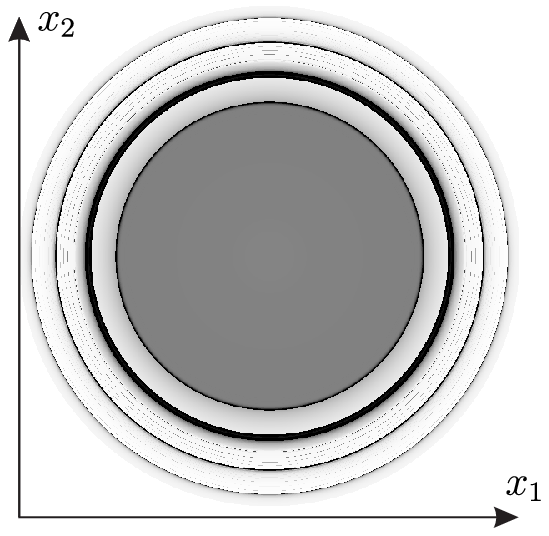}
        \caption{Bacteria density at end of experiment.}
        \label{figPetriPlate}
        \end{minipage}
\hspace{1.3cm}
\begin{minipage}[b]{0.4\textwidth}
        \centering
      \includegraphics[width=5cm]{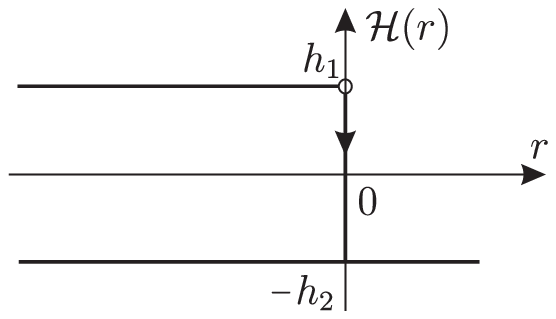}
        \caption{Hysteresis operator of relay type with threshold $0$.}
        \label{figRelay}
\end{minipage}
\end{center}
\end{figure}

In~\cite{Jaeger1,Jaeger2}, Hoppensteadt and J\"ager suggested a model for the growth of a colony of bacteria in a Petri plate. Bacteria fixed to
an agar gel in a Petri plate grow according to a hysteresis law. The growth rate depends on the relative
levels of the two diffusing substances: nutrient and by-product. A distinctive pattern of concentric rings
forms after a fixed amount of nutrient is added to the center of the plate, see Fig.~\ref{figPetriPlate}.
Mathematically, the Hoppensteadt--J\"ager  model is a system of reaction-diffusion equations coupled with an ordinary differential equation, and hysteresis enters the reaction terms. Numerical simulations carried
out in~\cite{Jaeger1,Jaeger2} were in accordance with the experiment. In a more general context, such systems may describe various hysteretic interactions between different diffusive and nondiffusive substances, see, e.g.~\cite{Marciniak06} for an application in the developmental biology.

Rigorous analysis of such systems appears to be
very nontrivial because the hysteresis may switch at different spatial points at different time moments.
Since 1980s, a series of results have been obtained about the existence of solutions, see~\cite{VisintinSpatHyst,Alt,VisintinBook,KopfovaSpatial1,AikiKopfova,Ilin,VisintinTenIssues14}. However,
questions about the uniqueness of solutions and their continuous dependence on initial data as well as
a thorough analysis of pattern formation remained open. For the purposes of this paper, we consider a simplified model given by the following scalar reaction-diffusion equation on a one-dimensional domain:
\begin{align}
& v_\tau=v_{xx}+\cH(v),\quad \tau>0,\ x\in(-1,1), \label{eqRDE}\\
& v|_{\tau=0}=\varphi(x),\quad x\in(-1,1),\label{eqRDEInitCond}
\end{align}
where $v=v(x,\tau)$ and $\cH(\cdot)$ is the hysteresis operator of relay type, see Fig.~\ref{figRelay}. To define it, we fix $h_1>0$ and $h_2\ge 0$ and set for any function $r\in C[0,\infty)$
\begin{equation}\label{nonideal}
 \cH(r)(\tau) =
 \begin{cases} h_1 & \text{if } r(s)< 0 \ \text{for all } s\in [0,\tau], \\
              -h_2 & \text{if } r(s)\ge 0 \ \text{for some}\  s\in
              [0,\tau].
 \end{cases}
\end{equation}
Note that the hysteresis operator is {\em rate-independent}, i.e., for any continuous increasing $g:[0,\infty)\to[0,\infty)$ with $g(0)=0$, we have
\begin{equation}\label{eqRateIndependent}
\cH(r(g(\cdot)))(\tau)=\cH(r)(g(\tau)),\quad \tau\ge 0.
\end{equation}

For a function $v(x,\tau)$ that is continuous with respect to $\tau$ for each fixed $x$, we set $\cH(v)(x,\tau):=\cH(v(x,\cdot))(\tau)$ for each fixed $x$. For a function $\varphi(x)$, we also set $\cH(\varphi)(x):=h_1$ if $\varphi(x)<0$ and $\cH(\varphi)(x):=-h_2$ if $\varphi(x)\ge 0$.

Even in such a simple form, problem~\eqref{eqRDE}, \eqref{eqRDEInitCond} reproduces the main difficulties of the original model. It turns out that the dynamics of $v$ and $\cH(v)$ essentially depends on the so-called (spatial) transversality condition. We say that a function $\phi(x)$ is {\it transverse with respect to $\cH(\varphi)$} provided the following holds: if $\varphi(x_0)=\varphi'(x_0)=0$ for some $x_0\in[-1,1]$, then $\cH(\varphi)(x_0)=-h_2$ in a neighborhood of $x_0$. We say that $v(x,\tau)$ is {\it transverse with respect to $\cH(v)$ on a time interval $I$} if $v(\cdot,\tau)$ is transverse with respect to $\cH(v)(\cdot,\tau)$ for each $\tau\in I$.
It was proved in~\cite{GurTikhNonlinAnal12,GurTikhSIAM13} that, for any transverse $\varphi(x)$, the solution $v(x,\tau)$ of problem~\eqref{eqRDE}, \eqref{eqRDEInitCond} exists, is transverse, is unique, and continuously depends on $\varphi$ on a time interval $(0,T)$. Moreover, the solution can be uniquely extended to its maximal interval of transverse existence $(0,T_{\rm max})$, where  either $T_{\rm max}=\infty$ or $T_{\rm max}$ is finite and the solution becomes nontransverse for $t=T_{\rm max}$. In~\cite{GurTikhMathBohem14}, these results were generalized to systems of reaction-diffusion equations (including the one suggested by Hoppensteadt and J\"ager). The case of a multi-dimensional domain was addressed in~\cite{CurranMastersThesis14}, and regularity of solutions was considered in~\cite{ApuUra-2015-regularity,ApuUra-2015-regularity2}, see also the survey~\cite{CurranGurTikhSurvey16}.

Although generically the initial data is transverse, it may become nontransverse in a finite time. In particular, this is the case for the Hoppensteadt--J\"ager model, where numerics indicates that all the concentric rings at Fig.~\ref{figPetriPlate} get formed via a nontransversality. To understand problem~\eqref{eqRDE}, \eqref{eqRDEInitCond} in the nontransverse case, we fix some $c>0$ and take $\varphi(x)=-cx^2+o(x^2)$ in a neighborhood $\mathcal B(0)$ of
$x=0$. By choosing a smaller neighborhood if needed, we have   $\varphi(x)<0$ for $x\in \mathcal B(0)\setminus\{0\}$. Now  for any $\ep > 0$, setting
$v_n(\tau;\varepsilon):=v(\varepsilon n,\tau)$, we replace the
continuous model~\eqref{eqRDE}, \eqref{eqRDEInitCond} in $\mathcal B(0)$ by the discrete one
\begin{equation*}%\label{eqDiscretePrototype}
\left\{
\begin{aligned}
& \dfrac{d v_n}{d \tau}  = \dfrac{\Delta v_n}{ \varepsilon^2}+\cH(v_n),\quad \tau>0, \ n=-N_\varepsilon,\dots,N_\varepsilon,\\
& v_n(0)=  - c(\varepsilon n)^2 + o(\varepsilon^2 n^2), \quad
n=-N_\varepsilon,\dots,N_\varepsilon,
\end{aligned}\right.
\end{equation*}
 where $\Delta v_n:=v_{n-1}-2v_n+v_{n+1}$ and $N_\varepsilon\to\infty$ as
$\varepsilon\to 0$. Since we are interested in small $\varepsilon$ and in the behavior in a small neighborhood $\mathcal B(0)$, we consider the
next approximation by omitting $o(\varepsilon^2 n^2)$ in the initial data and replacing $N_\varepsilon$ by $\infty$.
This yields
\begin{equation}\label{eqDiscretePrototypeZ}
\left\{
\begin{aligned}
& \dfrac{d v_n}{d \tau}  = \dfrac{\Delta v_n}{ \varepsilon^2}+\cH(v_n),\quad \tau>0,\ n\in\bbZ,\\
& v_n(0)=  - c(\varepsilon n)^2, \quad n\in\bbZ.
\end{aligned}\right.
\end{equation}
By our definition of the relay $\cH$, we have $\cH(v_0)(\tau)=-h_2$ for all $\tau\ge 0$, while, for $n\ne 0$, $\cH(v_n)(\tau)=h_1$ unless the input $v_n(\tau)$
achieves the zero threshold; if this happens at a finite time moment, the relay  {\it
switches} at this moment and since then the relay output $\cH(v_n)(\tau)$ equals $-h_2$.

\begin{figure}[ht]
%{\ \hfill\epsfxsize100mm\epsfbox{rattling.eps}\hfill\ }
\centering
\includegraphics[width=0.75\linewidth]{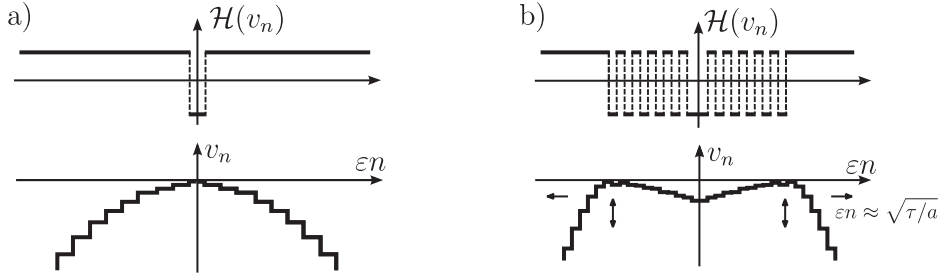}
\caption{The upper graphs represent spatial profiles of the hysteresis
$\cH(v_n)$ and the lower graphs the spatial profiles of the solution
$v_n$. {\rm a)} Nontransverse initial data. {\rm b)} Spatial profiles at a moment
$\tau>0$ for $h_1=h_2 > 0$.}
\label{figRattling}
\end{figure}

A nontrivial dynamics occurs in the case  $h_1>2c>0\ge -h_2$. To
indicate the difficulty, note that
$$
\dfrac{d
v_0}{d\tau}(0;\varepsilon)=-h_2-2c<0,\qquad \dfrac{d
v_n}{d\tau}(0;\varepsilon)=h_1-2c>0\ \text{for } n\ne0.
$$
Thus, for small $\tau>0$, the {\em node} $v_0(\tau;\varepsilon)$ decreases, while
all the other {\em nodes} $v_n(\tau;\varepsilon)$,
$n\ne0$, increase. Now it is not clear which
node achieves the threshold~$0$ and switches first and hence
what a further dynamics is. On the other hand, numerics provides quite a specific behavior of the solution, illustrated by Fig.~\ref{figRattling}.
The profile given by $v_n$ forms two humps symmetric to each other with respect to the origin. The humps propagate away from the origin, while the whole profile moves downwards and upwards (never going beyond the threshold $0$). Whenever the profile touches the threshold $0$, the corresponding node switches.   If we denote by $p(n)$ and $q(n)$ the
number of nodes in the set $\{v_1,\dots,v_{n}\}$ that
switch and do not switch, respectively, on the time interval
$[0,\infty)$, then numerics suggests that
\begin{equation}\label{eqDiscretePattern}
\lim\limits_{n\to\infty}\dfrac{q(n)}{p(n)}=\dfrac{h_2}{h_1}.
\end{equation}
Moreover,
\begin{equation}\label{eqPeriodicPattern}
\begin{aligned}
& \text{if $h_2/h_1=p_2/p_1$, where $p_1$ and $p_2$ are
co-prime integers, then}\\
&\text{for any $j$ large enough, the set
$\{v_{j+1},\dots, v_{j+p_1+p_2}\}$ contains exactly $p_1$ nodes}\\
& \text{that switch and $p_2$ nodes that do not switch on the time
interval $[0,\infty)$.}
\end{aligned}
\end{equation}
The switching pattern satisfying~\eqref{eqDiscretePattern} or, more specifically,~\eqref{eqPeriodicPattern} is called a {\em rattling pattern} or simply {\em rattling}.

We emphasize that the properties~\eqref{eqDiscretePattern} and~\eqref{eqPeriodicPattern} do not depend on $\varepsilon$.
In fact, in Section~\ref{secExistUnique}, we will show that $\varepsilon$
in~\eqref{eqDiscretePrototypeZ} can be scaled out, see~\eqref{eqRescale}. In particular, all the numerical
observations concerning the dynamics of $v_n$ have been done for
$\varepsilon=1$ and then transferred to an arbitrary $\varepsilon$
according to the scaling in~\eqref{eqRescale}. Therefore, the rattling describes a microstructure visible on the spatial scale of order $\varepsilon$, see the upper graph at Fig.~\ref{figRattling}.b. A similar effect is observable on different $d$-dimensional lattices, $d\ge 2$, and in slow-fast systems.
We comment on these issues and discuss a continuous limit as $\varepsilon\to 0$ in~\cite[Section~1]{GurTikhRattling1-0} and~\cite{CurranGurTikhSurvey16}. In connection with passing to the limit, let us also mention~\cite{HelmersHerrmannMMS13}, where a hysteretic behavior of forward-backward diffusion equations is analyzed via the spatial discretization.

% We also postpone a discussion of a continuous limit as $\varepsilon\to 0$ until Section~\ref{secDiscussionConclusion}.

The main goal of this paper is to analyse the property~\eqref{eqDiscretePattern} and to determine the asymptotic speed of propagation of the above mentioned humps formed by the solution profile, see Fig.~\ref{figRattling}.b. Let $\tau_n=\tau_n(\varepsilon)$ be the switching moment
of the node $v_n(\tau;\varepsilon)$ if this node switches  on the
time interval $[0,\infty)$. Then numerics shows that,
for any fixed $h_2 \ge 0$, the $\tau_n$'s that correspond to the switching nodes satisfy,
as $n\to\infty,$
\begin{equation}\label{eqTauAsymp}
\tau_n=a(\varepsilon n)^2+
\begin{cases}
\varepsilon^{2}
O(\sqrt{  n}) & \text{if } h_2=0,\\
\varepsilon^2 O(  n) & \text{if } h_2 > 0,
\end{cases}
 \end{equation}
where  $a>0$ and $O(\cdot)$ does not depend on~$\varepsilon$.

The case $h_2=0$ was studied in detail in~\cite{GurTikhRattling1-0}, where it was proved that all the nodes switch (which is certainly consistent with~\eqref{eqDiscretePattern}) and the constant $a$ from~\eqref{eqTauAsymp} was found. The case $h_2>0$ appears much more complicated, especially if $h_2/h_1$ is irrational. In this paper, using asymptotic methods developed in~\cite{GurTikhRattling1-0}, we find sufficient conditions that guarantee~\eqref{eqDiscretePattern}. More specifically, for any $h_2>0$, we will prove that if the switching nodes are {\em quasi-uniform}  (Condition~\ref{condPattern}) and  the switching moments grow asymptotically quadratically as in~\eqref{eqTauAsymp} (Condition~\ref{condAsymptk}), then the property~\eqref{eqDiscretePattern} takes place and the constant $a$ coincides with that found in~\cite{GurTikhRattling1-0} for $h_2=0$. In particular, we will show that $a$
depends on $h_1/c$ but {\em does not depend} on $h_2$ or $\varepsilon$. This is the contents of Theorem~\ref{thMain}. We believe that its proof contains essential ingredients that will allow us to prove the existence of the proper switching pattern in the spirit of~\cite{GurTikhRattling1-0}, i.e., without a priori assumptions about infinitely many switchings. This is a subject of future work.

From the viewpoint of applications, Theorem~\ref{thMain} may serve as a parameter identification tool. For instance, if one knows (or prepares ad hoc) the tangency constant $c$ and, by observing the emerging pattern, finds the constant $a$ in~\eqref{eqTauAsymp} and the frequency with which the switching nodes enter the pattern, then one can determine the parameters $h_1$ and $h_2$ of the hysteresis, using Theorem~\ref{thMain}.

The paper is organized as follows. In Section~\ref{secExistUnique}, we reduce problem~\eqref{eqDiscretePrototypeZ} to the problem with $\varepsilon=1$ by rescaling time $t :=\varepsilon^{-2} \tau$ and the solution itself $u_n(t):=\varepsilon^{-2}
v_n(\tau;\varepsilon)$. Then the switching moments $\tau_n$ (see~\eqref{eqTauAsymp}) of $v_n(\tau;\varepsilon)$ get transformed to the switching moments
\begin{equation}\label{eqtAsymp}
t_n=a n^2+
\begin{cases}

O(\sqrt{  n}) & \text{if } h_2=0,\\
 O(  n) & \text{if } h_2 > 0
\end{cases}
\end{equation}
of $u_n(t)$. Next, we prove the existence and uniqueness of the solution $u_n(t)$ for the rescaled problem (which does not depend on $\varepsilon$ any more). In Section~\ref{secMainResult}, we specify what {\em quasi-uniform switchings} mean (Condition~\ref{condPattern}) and formulate the main theorem (Theorem~\ref{thMain}) concerning the property~\eqref{eqDiscretePattern} and the constant $a$ in~\eqref{eqTauAsymp} and~\eqref{eqtAsymp}. In Section~\ref{secExample}, we discuss the notion of quasi-uniform switchings in more detail and give examples of periodic patterns that are observed numerically for rational $h_2/h_1$, and their generalizations, which we call quasiperiodic. In both cases, we show that the switching nodes are quasi-uniform (i.e., Condition~\ref{condPattern} holds and our main theorem applies). The rest of the paper is devoted to the proof of our main theorem. In Section~\ref{secElemntaryEstimates}, we prove preliminary estimates for
$u_n$, $\dot u_n$, $\Delta u_n:=u_{n-1}-2u_n+u_{n+1}$, and $\nabla u_n:=u_{n+1}-u_n$. In particular, we show the boundedness of $|\nabla u_n(t_n)|$ uniformly with respect to $n\in\bbZ$. Sections~\ref{secAsympUn} and~\ref{secAsympUnGrad} are central in the proof of our main result. In Section~\ref{secAsympUn}, we prove that a limit $p_*:=\lim\limits_{n\to\infty}p(n)/n$ exists and obtain an asymptotic expansion of $u_n(t_n)$ as $n\to\infty$. This expansion is based on the series representation of the solution $u_n(t)$ via the discrete Green function $\Gr_n(t)$, the asymptotic properties of the latter, and a proximity of emerging series to certain integrals. In Section~\ref{secAsympUnGrad}, we obtain an asymptotic expansion of $\nabla u_n(t_n)$ in a similar way.
Since $u_n(t_n)=0$ and $|\nabla u_n(t_n)|$ is bounded, the leading order terms in the asymptotic expansions of $u_n(t_n)$ and $\nabla u_n(t_n)$ must vanish. This yields two equalities for $a$ and $p_*$. Finally, in Section~\ref{secProofMainResult}, we show that these equalities imply that $p_*=h_1/(h_1+h_2)$ (which is equivalent to~\eqref{eqDiscretePattern}) and provide the unique possibility for choosing~$a$ (which appears to be independent of $h_2$). This concludes the prove of the main theorem. The paper also contains Appendix~\ref{secAppendixAsymptoticsFS}, where   discrete Green function $\Gr_n(t)$ is defined and some facts concerning its asymptotics are collected. 

%% file: sect2.tex
\section{Existence and uniqueness of solutions}\label{secExistUnique}

For a sequence $\{w_n\}_{n\in\bbZ}$ of real numbers, we will use throughout the
notation
\begin{equation}\label{eqNotationGrad}
\nabla w_n:=w_{n+1}-w_n,\qquad \Delta w_n:=\nabla w_{n}-\nabla
w_{n-1}=w_{n-1}-2w_n+w_{n+1}.
\end{equation}

We begin with scaling out $\varepsilon$ in problem~\eqref{eqDiscretePrototypeZ}. Setting
\begin{equation}\label{eqRescale}
 t :=
\varepsilon^{-2} \tau,\quad   u_n(t):=\varepsilon^{-2}
v_n(\tau;\varepsilon)
\end{equation}
and using the rate-independence of the hysteresis (see~\eqref{eqRateIndependent}), we obtain
$$
\cH(v_n)(\tau) = \cH(\varepsilon^2
u_n(\varepsilon^{-2}\cdot))(\tau) =
\cH(u_n(\varepsilon^{-2}\cdot))(\tau) =
\cH(u_n)(\varepsilon^{-2}\tau) = \cH(u_n)(t).
$$
Hence, problem~\eqref{eqDiscretePrototypeZ} can be rewritten as follows:
\begin{align}
& \dot u_n=\Delta u_n+\cH(u_n),\quad t>0,\ n\in\bbZ, \label{equ_nEquation1}\\
& u_n(0)=-c n^2,\quad n\in\mathbb Z, \label{equ_nEquation2}
\end{align}
where $\dot{}=d/d t$. Problem \eqref{equ_nEquation1}, \eqref{equ_nEquation2} does not involve $\ep$,
which justifies the fact that $u_n(t)$
in~\eqref{eqRescale} does not depend on $\varepsilon$. 

%Note that $c$
%in~\eqref{equ_nEquation2} could be also scaled out by replacing $u_n(t)$, $h_1$ and $-h_2$ by $c\tilde u_n(t)$, $c\tilde h_1$ and $-c\tilde h_2$, respectively.
%We prefer not to do this, in order to keep track of what
%exactly is influenced in our intermediate calculations by the
%``tangency'' constant~$c$.

We assume  that the following condition holds.

\begin{condition}\label{condHc}
$-h_2\le 0<2c<h_1$.
\end{condition}

We note that the function $\cH(r)(t)$ may have discontinuity
(actually, at most one) even if $r\in C^\infty[0,\infty)$.
Therefore, one cannot expect that a solution of
problem~\eqref{equ_nEquation1}, \eqref{equ_nEquation2} is
continuously differentiable on $[0,\infty)$. Thus, we define a
solution as follows.

\begin{definition}\label{defSolution}
  We say that a sequence
$\{u_n(t)\}_{n\in\bbZ}$ is a {\em solution of
problem~\eqref{equ_nEquation1}, \eqref{equ_nEquation2} on the time
interval $(0,T)$}, $T>0$,  if
\begin{enumerate}
\item $u_n\in C[0,T]$ for all $n\in\bbZ$,

\item\label{defSolution2} for each $t\in[0,T]$, there are
constants $A,\alpha\ge 0$ such that
$\sup\limits_{s\in[0,t]}|u_n(s)|\le A e^{\alpha |n|}$ for all
$n\in\bbZ$,

\item there is a finite  sequence $
0=\theta_0<\theta_1<\dots<\theta_J=T$, $J\ge1$,  such that $u_n\in
C^1[\theta_j,\theta_{j+1}]$ for all $n\in\bbZ$ and $j=0,\dots,J-1$,

\item the equations in~\eqref{equ_nEquation1} hold in
$(\theta_j,\theta_{j+1})$ for all $n\in\bbZ$ and  $j=0,\dots,J-1$,

\item $u_n(0)=-cn^2$ for all $n\in\bbZ$.
\end{enumerate}
  We say that a sequence
$\{u_n(t)\}_{n\in\bbZ}$ is a {\em solution of
problem~\eqref{equ_nEquation1}, \eqref{equ_nEquation2} on the time
interval $(0,\infty)$} if it is a solution on $(0,T)$ for all
$T>0$.
\end{definition}

\begin{remark}
If $\{u_n(t)\}_{n\in\bbZ}$ is a solution, then, as we have
mentioned above, the function $\cH(u_n)(t)$ has at most one
discontinuity point for each fixed $n\in\bbZ$. Hence, the
equations in~\eqref{equ_nEquation1} imply that each function $\dot
u_n(t)$ has at most one discontinuity point on $[0,\infty)$.
\end{remark}

For the
solution, we use the phrasing a {\em node $n$ switches $($at a moment $t$$)$} if $\cH(u_n)$ is discontinuous (at this moment $t$). If a node switches on the time interval $[0,\infty)$, we   call it a {\em switching node}; otherwise we call it a {\em nonswitching node}. For a switching node $n$, we denote the time moment at which it switches by $t_n$. By $\cS(t)$ we denote the set of nodes that   switch on the time
interval $[0,t]$, i.e.,
\begin{equation}\label{eqSt}
\cS(t):=\{k\in\bbZ: \cH(u_k)(t)=-h_2\},
\end{equation}
and let $\#(\cS(t))$ be the number of elements in $\cS(t)$.

The next result provides the existence and uniqueness of solutions as well as a two-sided bound for the switching moments $t_n$. Its proof uses the discrete Green function $\Gr_n(t)$. We refer the reader to Appendix~\ref{secAppendixAsymptoticsFS} for its definition and asymptotic properties, which play a central role in this paper.

\begin{theorem}\label{thExistenceUniqueness} Let Condition~$\ref{condHc}$ hold. Then the following is true.
\begin{enumerate}
 \item\label{thExistenceUniqueness1} Problem~\eqref{equ_nEquation1}, \eqref{equ_nEquation2}
has a unique solution $\{u_n(t)\}_{n\in\bbZ}$ on the time interval
$(0,\infty)$. It is given by
\begin{equation}\label{eqSolu}
u_n(t)=-cn^2+(h_1-2c)t-(h_1+h_2)\sum\limits_{k\in
\cS(t)}\Gr_{n-k}(t-t_k),\quad t\in[0,\infty),
\end{equation}
where $\Gr_n(t)$ is the discrete Green function given by~\eqref{eqynIntegralFormula}. The solution is symmetric$:$ $u_n(t)=u_{-n}(t)$ for all $n\in\bbZ$.

\item\label{thExistenceUniqueness2} For each switching node $n$, we have
$
t_n\ge \dfrac{c n^2}{h_1-2c}.
$

\item\label{thExistenceUniqueness3} The set $\cS(t)$
is finite for each $t\ge0$, symmetric with respect to the origin,
$\#(\cS(t))\to\infty$ as $t\to\infty$, and, for each switching node $n$, we have $t_{-n}=t_n$.
\end{enumerate}
\end{theorem}
\proof
Items~\ref{thExistenceUniqueness1} and~\ref{thExistenceUniqueness2} are proved in~\cite[Theorem~2.5]{GurTikhRattling1-0} in the case $h_2=0$. Their proof for $h_2>0$ is basically the same. Finiteness of $\cS(t)$ for each fixed $t\ge0$ follows from item~\ref{thExistenceUniqueness2}. Its symmetry and the relation $t_{-n}=t_n$ follow from the symmetry of the solution. Finally, assume that there are only finitely many, say $n_0\in\bbN$, switchings on the time interval $[0,\infty)$. Denote the last switching moment by $t_*\ge 0$.  Then  the solution formula~\eqref{eqSolu} contains the same $n_0$ terms in the sum for all $t>t_*$. Hence, since $h_1-2c>0$,  Proposition~\ref{theoremyn} implies that $u_n(t)\to\infty$ as $t\to \infty$ for all $n\in\bbZ$, which shows that the $(n_0+1)$th switching must happen. This proves that $\#(\cS(t))\to\infty$ as $t\to\infty$.
\endproof

%% file: sect3.tex
\section{Main result}\label{secMainResult}

First, we formulate assumptions meaning that the nodes switch according to a certain general pattern. Taking these assumptions for granted, we will rigorously determine specific parameters of this pattern. In particular, we will justify the conjecture~\eqref{eqDiscretePattern} and explicitly find the leading order term for the pattern propagation speed.

Denote by
$$
\cK=\{k_0,\pm k_1,\pm k_2, \pm k_3,\dots\}
$$
the set of nodes that switch on the time interval $[0,\infty)$. Let $0=k_0<k_1<k_2<\dots$ be the corresponding nonnegative nodes. For $n\in\bbN$, set
$$
\cK_n:=\cK\cap [-n,n],\quad \cK_n^+:=\cK\cap [1,n]
$$
and
\begin{equation}\label{eqpn}
p(n):=\#\left(\cK_{n-1}^+\right),\quad q(n):=n-p(n).
\end{equation}
  Thus, $p(n)$ and $q(n)$ are the number of switching and nonswitching nodes, respectively, from the spatial interval $[1,n-1]$. Note that the set $\cK$ is infinite and, thus, $p(n)\to\infty$ as $n\to\infty$ by Theorem~\ref{thExistenceUniqueness}.
We assume that the switching nodes are {\em quasi-uniform} in the following sense (see Section~\ref{secExample} below for a further discussion and examples of patterns with this property).
\begin{condition}[quasi-uniformity]\label{condPattern}
   $
  \sup\limits_{i=-p(n),\dots, p(n)}\left|\dfrac{k_i}{n}-\dfrac{i}{p(n)}\right|\to 0
  $ as $n\in\cK$, $n\to\infty$.
  \end{condition}

In what follows, we will use the notation
\begin{equation}\label{eqxynk}
x_{n,k_i}:=\dfrac{k_i}{n},\quad y_{n,i}:=\dfrac{i}{p(n)}
\end{equation}
for $n\in\cK$ and $i=-p(n),\dots,p(n)$.
Obviously, $y_{n,i}$ are uniformly distributed on the interval $[-1,1]$, while $x_{n,k_i}$ are generally not.
Condition~\ref{condPattern} is equivalent to the following:
\begin{equation}\label{eqCondPattern1Equiv}
\sup\limits_{i=-p(n),\dots, p(n)}\left|x_{n,k_i}-y_{n,i}\right|=0\quad\text{as } n\in\cK,\ n\to\infty,
\end{equation}
which motivates the notion of quasi-uniformity.

Our second assumption concerns the asymptotics of the switching moments $t_n$ and is also motivated by numerics, cf.~\eqref{eqTauAsymp}.

\begin{condition}\label{condAsymptk}
\begin{enumerate}
\item\label{condAsymptk1}
There exist  $a>0$ and $A>0$ such that the switching moment $t_k$ of each node $k\in\cK$ satisfies
\begin{equation}\label{eqEn+S1}
 t_k=a k^2+\omega_k,\quad |\omega_k|\le A |k|.
\end{equation}

\item\label{condAsymptk2} There exists $N_0\in\bbN$ and $\varepsilon_0>0$ such that, for any $n\in\cK\setminus[-N_0,N_0]$ and $k\in\cK_{n-1}$, we have
\begin{equation}\label{eqEn+S2}
 t_n-t_k\ge \varepsilon_0 (n^2-k^2).
\end{equation}
\end{enumerate}
\end{condition}

%\medskip
%
%{\bf REMARK FOR US. Suppose we assume $\omega_k=o(k^2)$ in~\eqref{eqEn+S1}. Then the asymptotics for $u_n(t_n)$ still works, but the one for $\nabla(u_n)(t_n)$ does not. }
%
%\medskip

Our main result guarantees that if the switching nodes are quasi-uniform and the switching moments grow asymptotically quadratic, then, actually, the proportion between the number of the switching  and nonswitching nodes satisfies conjecture~\eqref{eqDiscretePattern}. Furthermore, the coefficient $a$ in the leading order term of the asymptotics~\eqref{eqEn+S1} for~$t_k$ can be found explicitly.

\begin{theorem}\label{thMain}
  Let Conditions~$\ref{condHc}$, $\ref{condPattern}$, and~$\ref{condAsymptk}$ hold. Then the following hold.
 \begin{enumerate}
 \item There exist limits \begin{equation}\label{eqExistLimitpq}
  p_*:=\lim\limits_{n\to\infty}\dfrac{p(n)}{n},\quad q_*:=\lim\limits_{n\to\infty}\dfrac{q(n)}{n}.
  \end{equation}
 Moreover,  \begin{equation}\label{eqMain}
p_*=\dfrac{h_1}{h_1+h_2},\quad q_*=\dfrac{h_2}{h_1+h_2},\qquad \dfrac{q_*}{p_*}=\dfrac{h_2}{h_1}.
\end{equation}

 \item\label{thMain2} There exists an infinitely differentiable function $a_*:(2,\infty)\to (0,\infty)$ such that the coefficient $a$ in asymptotics~\eqref{eqEn+S1} is given by $a=a_*(h_1/c)$. In particular, $a$ does not depend on $h_2$. Furthermore, $a_*(\lambda)$ is monotone decreasing and satisfies the following relations$:$ \begin{equation}\label{AddTh}
     \mbox{$a_*(\lambda)\to +\infty$ as $\lambda\to 2$} \quad \mbox{and} \quad \mbox{$a_*(\lambda)\to 0$ as $\lambda\to +\infty$}.
     \end{equation}
\end{enumerate}
\end{theorem}

\begin{remark}
We emphasize  that item~\ref{thMain2} in Theorem~\ref{thMain} implies that the speed with which the pattern propagates does not depend on $h_2$, since $a=a_*(h_1/c)$ does not. The function $a_*(\cdot)$ is found in Section~\ref{secProofMainResult} as a unique root of each of the equations~\eqref{eqFa}, \eqref{eqGa}, \eqref{eqHa}.
Figure~\ref{figah1c} illustrates the dependence of $a_*$ on $h_1/c$.
\begin{figure}[ht]
\begin{center}
      \includegraphics[width=0.30\linewidth]{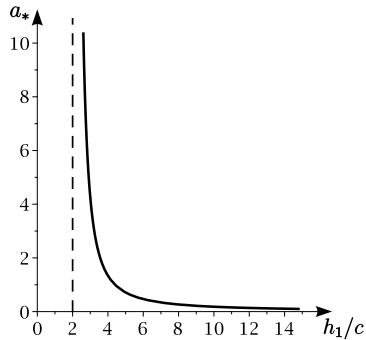}
\end{center}
\caption{Dependence of $a_*$ on $h_1/c$.} \label{figah1c}
\end{figure}
\end{remark}

The proof of Theorem~\ref{thMain} will be given in Section~\ref{secProofMainResult}, after some preparatory estimates in Section~\ref{secElemntaryEstimates} and asymptotic analysis in Sections~\ref{secAsympUn} and~\ref{secAsympUnGrad}.

 %based on asymptotic expansions of $u_n(t_n)$ and $\nabla u_n(t_n)$ and on a priori estimates of $\nabla u_n(t_n)$.  Section~\ref{secElemntaryEstimates} is devoted to a priori estimates of the solution $u_n(t)$ and, especially, of its gradient at the switching moment $t_n$. In sections~\ref{secAsympUn} and~\ref{secAsympUnGrad}, we obtain asymptotic expansions for $u_n(t_n)$ and $\nabla u_n(t_n)$, respectively. These formulas will be based on the properties of the discrete Green function $\Gr_n(t)$ (see Appendix~\ref{secAppendixAsymptoticsFS}). Finally, in section~\ref{secProofMainResult}, we combine all the previous results and prove Theorem~\ref{thMain}.

%% file: example.tex
\section{Discussion}\label{secExample}

\subsection{Discussion of Condition~$\ref{condPattern}$}

First, we observe the following. Suppose we are given an arbitrary set $\cK=\{k_0,\pm k_1,\pm k_2, \pm k_3,\dots\}$ such that $0=k_0<k_1<k_2<\dots$, not necessarily corresponding to the switching pattern for problem~\eqref{equ_nEquation1}, \eqref{equ_nEquation2}. Let us define the function $p(n)$ as in~\eqref{eqpn}. Then the following is true.

\begin{lemma}\label{lLimUniform}
  If  there exists a limit $\lim\limits_{n\in\cK, n \to \infty} {p(n)}/{n} > 0$, then
  \begin{equation*}%\label{P0}
     \sup\limits_{i=-p(n),\dots, p(n)} \left|\dfrac{k_i}{n} - \dfrac{i}{p(n)}\right| \to 0\quad\text{as } n\in\cK,\ n \to \infty.
  \end{equation*}
\end{lemma}
\proof
Obviously, it suffices to prove that
$\sup\limits_{j\in\cK_{n-1}^+}\left|\dfrac{j}{n}-\dfrac{p(j)}{p(n)}\right|\to 0$ as $n\in\cK$, $n\to\infty$. The latter follows by observing that
$$
\sup\limits_{j\in\cK_{n-1}^+}\left|\dfrac{j}{n}-\dfrac{p(j)}{p(n)}\right|= \sup\limits_{\alpha=j/n,\,j\in\cK_{n-1}^+}\left|\alpha\left(1-\dfrac{p(\alpha n)}{\alpha n}\dfrac{n}{p(n)}\right)\right|\to 0\quad\text{as }n\in\cK, n\to\infty.
$$
\endproof

Now we provide examples of patterns to which Theorem~\ref{thMain} applies. We consider {\em periodic} patterns that are observed numerically for rational $h_2/h_1$ and their generalizations, which we call {\em quasiperiodic}. In both cases, using Lemma~\ref{lLimUniform}, we show that Condition~\ref{condPattern} holds.

We fix   $\al \in (0, 1]$ and $\beta \in [0, 1)$. Denote by $\lfloor \cdot \rfloor$  the floor, or integer part, of a real number. The set
\begin{equation}\label{eqs1.1}
\cP:=\{0\} \cup \{n \in \ZZ\setminus \{0\}: \lfloor|n|\alpha + \beta\rfloor > \lfloor(|n|-1)\alpha + \beta\rfloor\}
\end{equation}
is said to be a {\it quasiperiodic  pattern}. If $\alpha$ is  rational, then we say that  $\cP$ is  a {\em periodic  pattern}.

\begin{figure}[ht]
  \centering
  \includegraphics[width=0.90\linewidth]{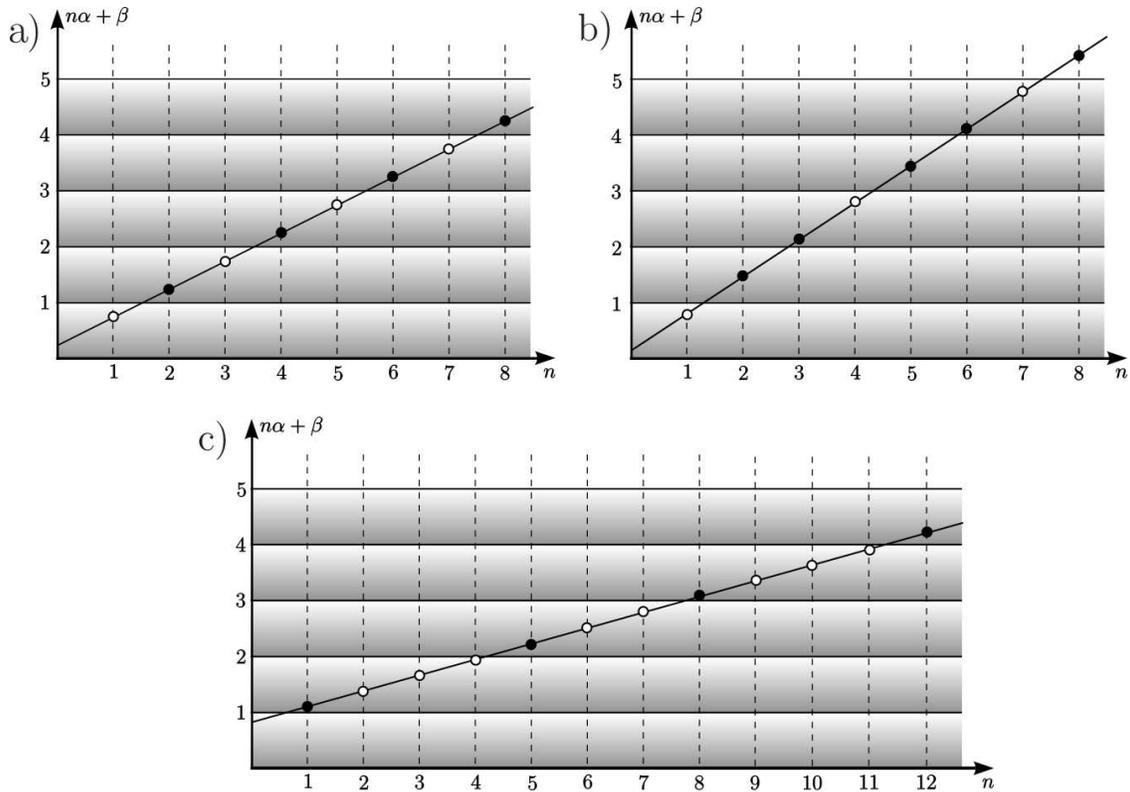}
  \caption{Periodic and quasiperiodic patterns$:$ $a)$ $\alpha = 1/2$, $\beta = 1/4;$ $b)$ $\alpha = 2/3$, $\beta = 1/5;$ $c)$~$\alpha = \sqrt{2}/5$, $\beta = 4/5$. The filled points of intersection correspond to the points included in~$\cP$}\label{figGrad}
\end{figure}

It is convenient to have the geometric description of the pattern $\cP$ in mind, see Fig.~\ref{figGrad}. On the plane $(x,y)$, consider the horizontal strips
\begin{equation}\label{eqStrips}
   \{0\le y<1\}, \quad \{1\le y<2\},\quad \{2\le y<3\},\dots
\end{equation}
and the points $Q_n=(n,n\alpha+\beta)$ lying on the line $y=\alpha x+\beta$. Then $0\in\cP$, and nodes $n>0$ and $-n<0$ belong to $\cP$  if and only if the points $Q_n$ and $Q_{n-1}$ do not belong to the same interval from~\eqref{eqStrips}.

Any switching pattern $\cK$ satisfying
$
\cK\setminus[-N_0,N_0]=\cP\setminus[-N_0,N_0]
$ for some $N_0\in\bbN$ necessarily consists of quasi-uniform nodes. Indeed,  $|p(n)-\lfloor\alpha n\rfloor|$ is bounded uniformly with respect to $n$. Hence, $\lim\limits_{n\to\infty}p(n)/n = \alpha$. Therefore, due to Lemma~\ref{lLimUniform},  the pattern $\cK$ satisfies Condition~\ref{condPattern}. Note that if  $\alpha$ is  rational, then  the co-prime integers $p_1$ and $p_2$ from the property~\eqref{eqPeriodicPattern} are defined by the equality $\alpha=p_1/(p_1+p_2)$.

\begin{remark}
It is easy to show that item~\ref{condAsymptk2} in Condition~\ref{condAsymptk} is a consequence of the following (simpler to formulate) condition: {\em There exists $N_0\in\bbN$ and $\varepsilon>0$ such that, for any $k_i\in \cK$, $k_i>N_0$, we have}
$ t_{k_i}-t_{k_{i-1}}\ge \varepsilon k_i.
$
\end{remark}

\subsection{Counter-example}

Now we show that the converse of Lemma~\ref{lLimUniform} is not true: quiasi-uniformity of the elements of a set $\cK$ (Condition~\ref{condPattern}) alone  does not guarantee the existence of a limit $\lim\limits_{n\to \infty}p(n)/n$. We will explicitly construct such sets~$\cK$. To show that they cannot be switching patterns generated by the solution of problem~\eqref{equ_nEquation1}, \eqref{equ_nEquation2}, an addition condition is necessary. In Theorem~\ref{thMain}, its role is played by Condition~\ref{condAsymptk}.

The idea of our construction is as follows. We   split the half-line $(0,\infty)$ into intervals $(M_j,M_{j+1}]$, $j=0,1,2,\dots,$ whose lengths $L_{j+1}:=M_{j+1}-M_j$ {\em grow fast} as $j\to\infty$. For each~$j$, we include in the set $\cK\cap (M_j,M_{j+1}]$ nodes with frequency $p_j\in(0,1)$. The frequencies~$p_j$ are chosen in such a way that they {\em slowly oscillate} between two different limits: $\liminf\limits_{j\to\infty}p_j$ and $\limsup\limits_{j\to\infty}p_j$. The whole construction is balanced in such a way that the nodes from $\cK$ are quasi-uniform (since the oscillation of $p_j$ is slow enough), but $p(n)/n$ has no limit (since $p_j$ has no limit).

Let us make the above argument rigorous. We consider a sequence of frequencies $p_j \in \bbQ$, $p_j \in (0, 1)$, possessing the {\em slow oscillation} property
\begin{equation}\label{eqFrequenciespj}
  \lim\limits_{j\to\infty}\dfrac{p_j}{p_{j+1}} = 1,\quad \liminf\limits_{j\to\infty} p_j < \limsup\limits_{j\to\infty} p_j
\end{equation}
and a sequence of interval lengths $L_j \in \bbN$ such that  $L_j p_j \in \bbN$ satisfying the {\em fast growth} property
\begin{equation}\label{eqLengths}
p_{j+1} L_{j+1} \geq (L_1 + \dots + L_j)(j+1).
\end{equation}
We define the intervals $(M_j,M_{j+1}]$ by setting $M_0:=0$ and $M_j := L_1 + \dots + L_j$. Finally, we define $\cK$ as the set consisting of $0$ and all the integer numbers  $n = M_j + r$ and $-n$ satisfying
\begin{equation}\label{2.3}
j \geq 0, \quad r \in [1, L_{j+1}], \quad \lfloor r p_{j+1} \rfloor > \lfloor (r-1) p_{j+1} \rfloor,
\end{equation}
cf.~\eqref{eqs1.1}.  In particular, for $R \in [1, L_{j+1}]$, we have
$
\#(\cK \cap (M_j, M_j + R]) = \lfloor Rp_{j+1} \rfloor,
$ which motivates the term ``frequency'' for $p_j$.

The following lemma asserts that the set $\cK$ possesses the desired properties.
\begin{lemma}
The set $\cK$ satisfies Condition~$\ref{condPattern}$, but
\begin{equation}\label{2.5.2}
\liminf\limits_{n\to\infty} \dfrac{p(n)}{n} \leq \liminf\limits_{j\to\infty} p_j, \quad \limsup\limits_{n\to\infty} \dfrac{p(n)}{n} \geq \limsup\limits_{j\to\infty} p_j.
\end{equation}
\end{lemma}
\proof 1. First, we prove that   $
  \sup\limits_{i=-p(n+1),\dots, p(n+1)}\left|\dfrac{k_i}{n}-\dfrac{i}{p(n+1)}\right|\to 0
  $ as $n\in\cK$, $n\to\infty$, which implies Condition~$\ref{condPattern}$ because $p(n)\to\infty$ as $n\to\infty$.
Set
$$
  l_j  := L_j p_j=\#(\cK \cap (M_{j-1}, M_j]), \quad m_j := l_1 + \dots + l_j=\#(\cK \cap (0, M_j]), % \quad m_0 = 0,
$$
$$
 b_j := \left|\dfrac{M_j p_{j+1}}{m_j} - 1\right|.
$$
It is easy to show that~\eqref{eqFrequenciespj} and~\eqref{eqLengths} yield
\begin{equation}\label{1.1}
  \lim\limits_{j\to\infty} \dfrac{M_j}{M_{j+1}} = 0, \quad  \lim\limits_{j\to\infty} \dfrac{m_j}{m_{j+1}} = 0, \quad  \lim\limits_{j\to\infty} \dfrac{p_jM_j}{m_j} = 1, \quad   \lim\limits_{j\to\infty} b_j = 0.
\end{equation}

Now, for $i=1,\dots, p(n+1)$, we represent
\begin{equation}\label{2.1}
\begin{aligned}
  n   &= M_j + r  \mbox{ with }  r \in [1, L_{j+1}], &  i  &= m_s + \rho  \mbox{ with }  \rho \in [1, l_{s+1}],\\
  p(n+1) & = m_j + \lfloor r p_{j+1} \rfloor, &  k_i  & = M_s + \left\lceil \dfrac{\rho}{p_{s+1}} \right\rceil,
  \end{aligned}
\end{equation}
where $\lceil x \rceil := \min_{u \in \bbZ}\{u \geq x\}$. In particular,
\begin{equation}\label{Add2.1}
s < j \quad \mbox{or} \quad (s = j \mbox{ and } \rho \leq r p_{j+1}).
\end{equation}

If $s \leq j-2$, then
\begin{equation}\label{4.1}
\dfrac{k_i}{n} \leq \dfrac{M_{s+1}}{M_j} \leq \dfrac{M_{j-1}}{M_j} \to 0\quad \text{as } {j \to \infty},
\end{equation}
\begin{equation}\label{4.2}
  \frac{i}{p(n+1)} \leq \dfrac{m_{s+1}}{m_j} \leq \dfrac{m_{j-1}}{m_j} \to  0\quad \text{as } {j \to \infty}.
\end{equation}

It remains to consider the case $s = j-1$ or $s = j$. We have
\begin{equation}\label{eqkin-ipn}
\frac{k_i}{n} - \frac{i}{p(n+1)} =
    \dfrac{M_s + \left\lceil \dfrac{\rho}{p_{s+1}} \right\rceil}{M_j + r} - \dfrac{m_s + \rho}{m_j + \lfloor rp_{j+1} \rfloor}.
\end{equation}
We rewrite the first term in the right-hand side in~\eqref{eqkin-ipn} as follows:
  \begin{equation}\label{3.1}
    \dfrac{M_s + \left\lceil \dfrac{\rho}{p_{s+1}} \right\rceil}{M_j + r} = \dfrac{M_s + \dfrac{\rho}{p_{s+1}}}{M_j + r}(1+ \delta_1),
  \end{equation}
  where
  \begin{equation}\label{3.2}
    |\delta_1| \leq \dfrac{1}{M_s}.
  \end{equation}
Similarly,
  \begin{equation}\label{3.3}
  \begin{aligned}
    \dfrac{M_s + \dfrac{\rho}{p_{s+1}}}{M_j + r} & =
    \dfrac{M_s p_{s+1}+ \rho}{M_jp_{s+1} + r p_{s+1}} =
    \dfrac{m_s \dfrac{M_s p_{s+1}}{m_s} + \rho}{\dfrac{p_{s+1}}{p_{j+1}}
    \left(m_j \dfrac{M_jp_{j+1}}{m_j} + rp_{j+1} \right)} \\
    & = \dfrac{p_{j+1}}{p_{s+1}} \cdot \dfrac{m_s + \rho}{m_j + rp_{j+1}} \cdot \dfrac{1+ \delta_2}{1+ \delta_3},
  \end{aligned}
\end{equation}
  where
  \begin{equation}\label{3.4}
  |\delta_2| \leq b_s, \quad |\delta_3| \leq b_j.
  \end{equation}
We also rewrite the second term in the right-hand side in~\eqref{eqkin-ipn} as follows:
  \begin{equation}\label{3.5}
    \dfrac{m_s + \rho}{m_j + \lfloor rp_{j+1} \rfloor} = \dfrac{m_s + \rho}{m_j + rp_{j+1}} \cdot \dfrac{1}{1+ \delta_4},
  \end{equation}
  where
  \begin{equation}\label{3.6}
    |\delta_4| \leq \dfrac{1}{m_j}.
  \end{equation}

Now, taking into account \eqref{Add2.1},  \eqref{eqkin-ipn}, \eqref{3.1}, \eqref{3.3}, and  \eqref{3.5}, we obtain
$$
\begin{aligned}
  \left|\frac{k_i}{n} - \frac{i}{p(n+1)}\right| & = \dfrac{m_s + \rho}{m_j + rp_{j+1}}\left|\dfrac{(1+\delta_1)(1+\delta_2)}{1+\delta_3} \cdot \dfrac{p_{j+1}}{p_{s+1}} - \dfrac{1}{1+\delta_4} \right| \\
  & \leq \left|\dfrac{(1+\delta_1)(1+\delta_2)}{1+\delta_3} \cdot \dfrac{p_{j+1}}{p_{s+1}} - \dfrac{1}{1+\delta_4} \right|.
\end{aligned}
$$
Since $s = j-1$ or $s = j$, relations~\eqref{eqFrequenciespj}, \eqref{1.1},  \eqref{3.2}, \eqref{3.4}, and \eqref{3.6} imply that the latter expression tends to 0 as $j \to \infty$.

2. In order to prove relations \eqref{2.5.2}, we consider the subsequence $n_j = M_j$. Note that due to \eqref{2.1} the equality ${p(n_j+1)}/{n_j} = {m_j}/{M_j}$ holds, and~\eqref{eqFrequenciespj} together with the third relation in~\eqref{1.1} imply~\eqref{2.5.2}.
\endproof

%% file: sect4.tex
\section{Elementary estimates}\label{secElemntaryEstimates}

We begin with some straightforward estimates for
$u_n$, $\dot u_n$, $\Delta u_n$, and $\nabla u_n$ (see notation~\eqref{eqNotationGrad}). In particular, our goal is to prove the following result about the boundedness of $\nabla u_n(t_n)$.

\begin{proposition}\label{lGradientuBounded}
Let Condition~$\ref{condHc}$ hold. Then, for each node $n\in\cK$, the following relations hold at the switching moment~$t_n$$:$
$$
\nabla u_n(t_n)\in[-(2h_1+h_2),0],\quad \nabla u_{n-1}(t_n)\in[0,(2h_1+h_2)].
$$
\end{proposition}

The proof will be given in the end of this section as a consequence of several lemmas (in which we assume that Condition~\ref{condHc} holds).

\begin{lemma}\label{luLessZero}
$u_n(t)\le 0$ for all $n\in\bbZ$, $t\ge 0$.
\end{lemma}
\proof Using the solution formula~\eqref{eqSolu} and the relations~\eqref{eqy_nGreen} for the discrete Green function $\Gr_n(t)$, it is easy to see that, on a sufficiently small time interval, all the nodes satisfy $u_m(t)\le 0$, $m\in\bbZ$. Let $\theta>0$ be a time moment at which at least one of the nodes, say $u_n$, is on the threshold $0$, while all the nodes still satisfy $u_m(\theta)\le 0$, $m\in\bbZ$. Then $\Delta u_n(\theta)\le 0$. If $\Delta u_n(\theta)<0$ or if  $\Delta u_n(\theta)=0$ and $-h_2<0$, then $\dot u_n(\theta+0)<0$ and $u_n(t)$ must decrease on $(\theta,\theta+\varepsilon)$ for a sufficiently small $\varepsilon>0$.

It remains to consider the case $h_2=0$,  $u_n(\theta)=0$ for some
$n\in\bbZ$, and $\Delta u_n(\theta)=0$. If $u_n(\theta)=0$ for all
$n\in\bbZ$, then Theorem~\ref{thExistenceUniqueness} implies that $u_n(t)\equiv 0$ for $n\in\bbZ$ and $t\in[\theta,\infty)$.

Assume that
$$
u_{n}(\theta)=u_{n\pm 1}(\theta)=u_{n\pm 2}(\theta)=\dots u_{n\pm
(j-1)}(\theta)=0
$$
for some $j\ge 2$, but  $u_{n+j}(\theta)<0$ or $u_{n-j}(\theta)<0$.
Then
$$\underbrace{\Delta\dots\Delta}_{s\
\text{times}}u_n(\theta)=0, \quad s=1,\dots,j-1;\qquad
\underbrace{\Delta\dots\Delta}_{j\ \text{times}}u_n(\theta)<0.
$$
Therefore, $\dfrac{d^su_n(\theta)}{dt^s}=0$ for $s=1,\dots,j-1$, but
$\dfrac{d^ju_n(\theta)}{dt^j}<0$ (where $d/dt$ is understood as the right
derivative).
\endproof

\begin{lemma}\label{ludotBounded}
$\dot u_n(t)\in[-(h_1+h_2),h_1-2c]$ for all $n\in\bbZ$, $t\ge 0$.
\end{lemma}
\proof 1. By Theorem~\ref{thExistenceUniqueness}, there is a sequence $0=:\theta_0 <\theta_1<\theta_2<\dots$ (of switching moments) such that $\lim\limits_{j\to\infty}\theta_j=\infty$ and there are no switchings for $t\in(\theta_j,\theta_{j+1})$, $j=0,1,2,\dots$.

Set $v_n(t):=\dot u_n(t)$ for $t>0$, $t\ne\theta_j$ ($j=0,1,2,\dots$). Note that $v_n\in C[\theta_j,\theta_{j+1}]$, $j=0,1,2,\dots$. However, if a node $n$  switches at some moment $\theta_j$, then $v_n(t)$ has a jump discontinuity at $\theta_j$.

At the initial moment, we
have
$$
\begin{aligned}
v_0(0)&=-2c-h_2\in [-(h_1+h_2),h_1-2c],\\
v_n(0)&=-2c+h_1\in[-(h_1+h_2),h_1-2c],\quad n\ne 0.
\end{aligned}
$$

Due to~\eqref{equ_nEquation1}, on the interval $[\theta_0,\theta_1)$,
$v_n(t)$ satisfies
$$
\dot v_n=\Delta v_n,\quad n\in\bbZ.
$$
Therefore, similarly to the proof of Lemma~\ref{luLessZero}, we
conclude that $v_n(t)\in [-(h_1+h_2),h_1-2c]$ for all $[\theta_0,\theta_1)$.

2. At the switching moment $\theta_1$, some nodes switch. Let $n$ be one of them. Then $v_n(\theta_1-0)\ge
0$ and $\Delta u_n(\theta_1)=\Delta u_n(\theta_1\pm0)\le 0$. Therefore,
$$
v_n(\theta_1+0)=\Delta u_n(\theta_1+0)-h_2\le -h_2< h_1-2c,
$$
$$
\begin{aligned}
 v_n(\theta_1+0)=\Delta u_n(\theta_1+0)-h_2&=\Delta
u_n(\theta_1-0)+h_1-(h_1+h_2)\\
&=v_n(\theta_1-0)-(h_1+h_2)\ge -(h_1+h_2).
\end{aligned}
$$

If the node $n$ does not switch at $\theta_1$, then $v_n(\theta_1+0)=v_n(\theta_1-0)\in [-(h_1+h_2),h_1-2c]$.

Now, as in part 1 of the proof, we conclude that $v_n(t)\in
[-(h_1+h_2),h_1-2c]$ for all $n\in\bbZ$ and $t\in [\theta_1,\theta_2)$. Arguing by induction, we
complete the proof.
\endproof

\begin{lemma}\label{lDeltauBounded}
$\Delta u_n(t)\in[-(2h_1+h_2),h_1+h_2-2c]$ for all $n\in\bbZ$, $t\ge
0$.
\end{lemma}
\proof The result follows from Lemma~\ref{ludotBounded} and
equation~\eqref{equ_nEquation1}.
\endproof

\proof[Proof of Proposition~$\ref{lGradientuBounded}$]
By Lemma~\ref{luLessZero}, at the switching moment  $t=t_n$ ($n\in\cK$), we have $\nabla u_n\le 0$
and $\nabla u_{n-1}\ge 0$. Therefore, using
Lemma~\ref{lDeltauBounded}, we have at the moment $t=t_n$
$$
\begin{gathered}
\nabla u_n=\nabla u_{n-1}+\Delta u_n\ge \Delta u_n\ge
-(2h_1+h_2),\\
\nabla u_{n-1}=\nabla u_n-\Delta u_n\le -\Delta u_n\le 2h_1+h_2.
\end{gathered}
$$
\endproof

%% file: sect5.tex
\section{Asymptotics for $u_n(t_n)$}\label{secAsympUn}

In this section and in the next one, we will use the standard notation $O(\cdot)$ and $o(\cdot)$. For example, given two sequences of real numbers $v_{n,k}$ and $w_{n,k}$  ($n,k\in\bbN$), we write $v_{n,k}=O(w_{n,k})$ if there exists a constant $C>0$ such that $|v_{n,k}|\le C |w_{n,k}|$ for all $n,k\in\bbN$. Given a sequence of real numbers $v_n$ ($n\in\bbN$), we write $v_n=o(n)$ if $v_n/n\to 0$ as $n\to\infty$.

The goal of this section is to prove the following result.

\begin{figure}[ht]
  \centering
  \includegraphics[width=0.75\linewidth]{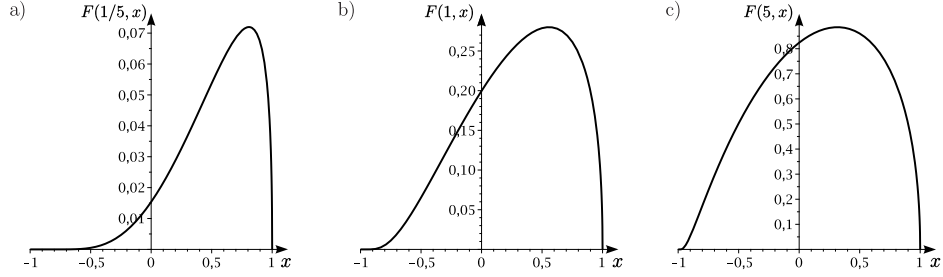}
  \caption{Graphs of function $F(a, \cdot)$ for $a)$ $a = 1/5$, $b)$ $a = 1$, $c)$ $a = 5$.}\label{figF}
\end{figure}

\begin{proposition}\label{propAsympun}
Let Conditions~$\ref{condHc}$, $\ref{condPattern}$, and~$\ref{condAsymptk}$ hold. Then there exists a limit
\begin{equation}\label{eqExistLimitpqProp}
  p_*:=\lim\limits_{n\to\infty}\dfrac{p(n)}{n},
  \end{equation}
  and the numbers  $a$ from Condition~$\ref{condAsymptk}$ and $p_*$ from~\eqref{eqExistLimitpqProp} satisfy the equality
\begin{equation}\label{eqa1}
-c+(h_1-2c)a-(h_1+h_2)p_* I_F(a)=0,
\end{equation}
where
\begin{equation}\label{eqF}
I_F(a):=\int\limits_{-1}^1 F(a,x)\,dx,\quad  F(a,x):=\sqrt{a}{\sqrt{1-x^2}}\,f\left(\dfrac{1}{\sqrt{a}}\sqrt{\dfrac{1-x}{1+x}}\right)
\end{equation}
$($see Fig.~$\ref{figF}$$)$, and the function $f$ is defined in~\eqref{eqfghft}.
\end{proposition}

\proof
For $n\in\cK$, $n>N_0$ ($N_0$ is the number from item~\ref{condAsymptk2} in Condition~\ref{condAsymptk}), we will rewrite the equality
$u_n(t_n)=0$. Using~\eqref{eqSolu}, we have
\begin{equation}\label{equntn1}
0=u_n(t_n)  =-cn^2+(h_1-2c)t_n-(h_1+h_2)\sum\limits_{k\in\cK_{n-1}} \Gr_{n-k}(t_n-t_k).
\end{equation}

Consider the sum in~\eqref{equntn1}. We write $t_n-t_k=a(n^2-k^2)+\omega_n-\omega_k$ according to Condition~\ref{condAsymptk} and expand $\Gr_{n-k}$ in the Taylor series about $a(n^2-k^2)$. Taking into account Proposition~\ref{theoremyn} and  notation~\eqref{eqxynk} for $x_{n,k_i}$, we obtain
\begin{equation}\label{equntn1'}
\begin{aligned}
\sum\limits_{k\in\cK_{n-1}}\Gr_{n-k}(t_n-t_k) & =\sum\limits_{k\in\cK_{n-1}}\Gr_{n-k}(a(n^2-k^2)) + \Sigma_{n}\\
& = n^2\sum\limits_{i=-p(n)}^{p(n)}\dfrac{1}{n} F(a,x_{n,k_i})
+ \sum\limits_{k\in\cK_{n-1}} O\left((n^2-k^2)^{-1/2}\right) +  \Sigma_{n},
\end{aligned}
\end{equation}
where  the function $F$ is defined in~\eqref{eqF},
\begin{equation}\label{eqSigma1n1}
\Sigma_{n} :=\sum\limits_{k\in\cK_{n-1}}\dot
\Gr_{n-k}(\sigma_{n,k})(\omega_n-\omega_k),
\end{equation}
and $\sigma_{n,k}\ge\min\{\varepsilon_0,a\} (n^2-k^2)$ due to~\eqref{eqEn+S2}. Note that
\begin{equation}\label{eqSumOn2k2}
\sum\limits_{k\in\cK_{n-1}} O\left((n^2-k^2)^{-1/2}\right)=O(1).
\end{equation}
Further, since $\omega_n-\omega_k=O(n)$ by~\eqref{eqEn+S1}, we obtain from~\eqref{eqSigma1n1} and Proposition~\ref{theoremyn}
\begin{equation}\label{equntn1''}
\begin{aligned}
\Sigma_{n}=\sum\limits_{k\in\cK_{n-1}}
O((n^2-k^2)^{-1/2})O(n)=O(n).
\end{aligned}
\end{equation}

Combining~\eqref{equntn1'}, \eqref{eqSumOn2k2}, and \eqref{equntn1''}, we obtain
\begin{equation}\label{equntn2}
\begin{aligned}
&\sum\limits_{k\in\cK_{n-1}}\Gr_{n-k}(t_n-t_k) = n^2\sum\limits_{i=-p(n)}^{p(n)}\dfrac{1}{n} F(a,x_{n,k_i})
+ O(n).
\end{aligned}
\end{equation}
Now we represent the sum in~\eqref{equntn2} as follows (using notation~\eqref{eqxynk} for $y_{n,i}$ and recalling that $F(-1)=0$):
\begin{equation}\label{equntn2'}
\begin{aligned}
&\sum\limits_{i=-p(n)}^{p(n)}\dfrac{1}{n} F(a,x_{n,k_i})\\
&\qquad =\dfrac{p(n)}{n}\left(\sum\limits_{i=-p(n)+1}^{p(n)}\dfrac{1}{p(n)} F(a,y_{n,i}) + \sum\limits_{i=-p(n)}^{p(n)}\dfrac{1}{p(n)}\left(F(a,x_{n,k_i})-F(a,y_{n,i})\right) \right).
\end{aligned}
\end{equation}
 The first sum in~\eqref{equntn2'} is the Riemann sum for the integral $I_F(a)$ defined in~\eqref{eqF}. Since $p(n)\to\infty$ as $n\to\infty$, we have
\begin{equation}\label{equntn2'Intermediate1}
\sum\limits_{i=-p(n)+1}^{p(n)}\dfrac{1}{p(n)} F(a,y_{n,i})=I_F(a)+o(1).
\end{equation}
Finally, taking into account Condition~\ref{condPattern} and the uniform continuity of $F(a,x)$ for $x\in[-1,1]$ (and a fixed $a$), we estimate the last sum  in~\eqref{equntn2'}  as follows:
\begin{equation}\label{equntn2'Intermediate2}
\begin{aligned}
&\left|\sum\limits_{i=-p(n)}^{p(n)}\dfrac{1}{p(n)}\left(F(a,x_{n,k_i})-F(a,y_{n,i})\right)\right|\\
&\qquad \le \dfrac{2p(n)+1}{p(n)} \sup\limits_{i=-p(n),\dots,p(n)}|F(a,x_{n,k_i})-F(a,y_{n,i})| = o(1).
\end{aligned}
\end{equation}
Due to~\eqref{equntn2'}--\eqref{equntn2'Intermediate2}, relation \eqref{equntn2} takes the form
\begin{equation}\label{equntn2''}
\sum\limits_{k\in\cK_{n-1}}\Gr_{n-k}(t_n-t_k) = n^2 \dfrac{p(n)}{n} I_F(a) + o(n^2).
\end{equation}
Combining relations~\eqref{equntn1} and~\eqref{equntn2''}, we conclude that
$$
\left(-c+(h_1-2c)a-(h_1+h_2) \dfrac{p(n)}{n} I_F(a)\right)n^2 + o(n^2)=0.
$$
Since $I_F(a)>0$, this implies that the limit $p_*:=\lim\limits_{n\in\cK,\, n\to\infty}\dfrac{p(n)}{n}$ exists and equality~\eqref{eqa1} holds.

To complete the proof, we have to show that the limit in~\eqref{eqExistLimitpqProp} exists along any sequence. Fix an arbitrary $\varepsilon>0$. Then there exists $N\ge 1/\varepsilon$ such that $\left|\dfrac{p(k_i)}{k_i}-p_*\right|\le \varepsilon$ for all $k_i\in\cK$, $k_i\ge N$. Now any $m\ge N$ belongs to some interval $(k_{i-1},k_i]$, and we have
$$
\begin{aligned}
\left|\dfrac{p(m)}{m}-p_*\right| & =\left|\dfrac{p(k_i)}{m}-p_*\right|\le \left|\dfrac{p(k_i)}{k_i}-p_*\right|+\left|\dfrac{p(k_i)}{k_i}-\dfrac{p(k_i)}{m}\right|\\
&\le \varepsilon + \dfrac{p(k_i)}{m} \left(1-\dfrac{m}{k_i}\right) \le \varepsilon + \dfrac{k_i-1}{k_{i-1}} \left(1-\dfrac{k_{i-1}}{k_i}\right).
\end{aligned}
$$
It remains to note that $\dfrac{k_{i}}{k_{i-1}}\to 1$ as $k_i\in\cK$, $k_i\to\infty$ due to Condition~\ref{condPattern}.
\endproof

%% file: sect6.tex
\section{Asymptotics for $\nabla u_n(t_n)$}\label{secAsympUnGrad}

In this section, we will prove the following result.

\begin{proposition}\label{propAsympGradun}
Let Conditions~$\ref{condHc}$, $\ref{condPattern}$, and~$\ref{condAsymptk}$ hold. Then the numbers  $a$ from Condition~$\ref{condAsymptk}$ and $p_*$ from~\eqref{eqExistLimitpqProp}   satisfy the equality
\begin{equation}\label{eqa2}
-2c-(h_1+h_2) p_* I_G(a)=0,
\end{equation}
where
\begin{equation}\label{eqG}
I_G(a):=\int\limits_{-1}^1 G(a,x)\,dx,\quad  G(a,x):=g\left(\dfrac{1}{\sqrt{a}}\sqrt{\dfrac{1-x}{1+x}}\right)
\end{equation}
$($see Fig.~$\ref{figG}$$)$, and the function $g$ is defined in~\eqref{eqfghft}.
\begin{figure}[ht]
  \centering
  \includegraphics[width=0.75\linewidth]{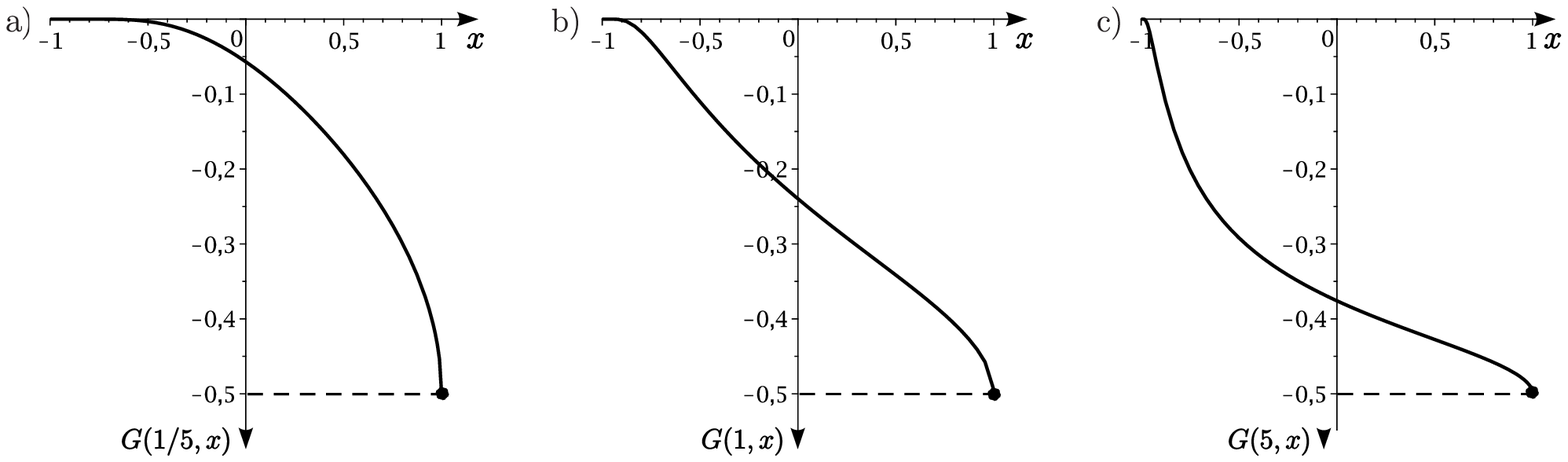}
  \caption{Graphs of function $G(a, \cdot)$ for $a)$ $a = 1/5$, $b)$ $a = 1$, $c)$ $a = 5$.}\label{figG}
\end{figure}
\end{proposition}
\proof
For $n\in\cK$, $n>N_0$ ($N_0$ is the number from item~\ref{condAsymptk2} in Condition~\ref{condAsymptk}), we will find an asymptotic expansion of $\nabla u_n(t_n)$. Using~\eqref{eqSolu}, we have
\begin{equation}\label{eqNablauntn1}
\begin{aligned}
\nabla u_n(t_n)& =-2cn-c-(h_1+h_2)\sum\limits_{k\in\cK_{n-1}}\nabla
\Gr_{n-k}(t_n-t_k).
\end{aligned}
\end{equation}

Consider the last sum in~\eqref{eqNablauntn1}. We write $t_n-t_k=a(n^2-k^2)+\omega_n-\omega_k$ according to Condition~\ref{condAsymptk} and expand $\nabla \Gr_{n-k}$ in the Taylor series about $a(n^2-k^2)$. Taking into account Proposition~\ref{theoremyn} and  notation~\eqref{eqxynk} for $x_{n,k_i}$, we obtain
\begin{equation}\label{equntn1'Grad}
\begin{aligned}
\sum\limits_{k\in\cK_{n-1}}\nabla \Gr_{n-k}(t_n-t_k) & =\sum\limits_{k\in\cK_{n-1}}\nabla \Gr_{n-k}(a(n^2-k^2)) + \Theta_{n}\\
&  = n\sum\limits_{i=-p(n)}^{p(n)}\dfrac{1}{n} G(a,x_{n,k_i})
+ \sum\limits_{k\in\cK_{n-1}} O\left((n^2-k^2)^{-1/2}\right) +  \Theta_{n},
\end{aligned}
\end{equation}
where the function  $G$ is defined in~\eqref{eqG},
\begin{equation}\label{eqSigma1n1Grad}
\Theta_{n} :=\sum\limits_{k\in\cK_{n-1}} \nabla\dot \Gr_{n-k}(\sigma_{n,k})(\omega_n-\omega_k),
\end{equation}
and $\sigma_{n,k}\ge\min\{\varepsilon_0,a\} (n^2-k^2)$ due to~\eqref{eqEn+S2}. Since $\omega_n-\omega_k=O(n)$ by~\eqref{eqEn+S1}, we obtain from~\eqref{eqSigma1n1Grad} and Proposition~\ref{theoremyn}
\begin{equation}\label{equntn1''Grad}
\begin{aligned}
\Theta_{n}=\sum\limits_{k\in\cK_{n-1}}
O((n^2-k^2)^{-1})O(n)=O(\ln n).
\end{aligned}
\end{equation}
Combining~\eqref{equntn1'Grad}, \eqref{eqSumOn2k2}, and \eqref{equntn1''Grad}, we obtain
\begin{equation}\label{equntn2Grad}
\begin{aligned}
&\sum\limits_{k\in\cK_{n-1}}\nabla \Gr_{n-k}(t_n-t_k) = n\sum\limits_{i=-p(n)}^{p(n)}\dfrac{1}{n} G(a,x_{n,k_i})
+ O(\ln n).
\end{aligned}
\end{equation}
Similarly to the proof of Proposition~\ref{propAsympun} (cf.~\eqref{equntn2'}--\eqref{equntn2'Intermediate2}), taking into account the existence of the limit in~\eqref{eqExistLimitpqProp}, we obtain
$$
\begin{aligned}
&\sum\limits_{i=-p(n)}^{p(n)}\dfrac{1}{n} G(a,x_{n,k_i}) \\
&\qquad = \dfrac{p(n)}{n}\left(\sum\limits_{i=-p(n)+1}^{p(n)}\dfrac{1}{p(n)} G(a,y_{n,i}) + \sum\limits_{i=-p(n)}^{p(n)}\dfrac{1}{p(n)}\left(G(a,x_{n,k_i})-G(a,y_{n,i})\right) \right)\\
& \qquad =p_* I_G(a)+o(1),
\end{aligned}
$$
where $p_*$ is given by~\eqref{eqExistLimitpqProp} and the integral $I_G(a)$ is defined in~\eqref{eqG}.
Therefore, \eqref{equntn2Grad} takes the form
\begin{equation}\label{equntn2''Grad}
\sum\limits_{k\in\cK_{n-1}}\nabla \Gr_{n-k}(t_n-t_k) = n p_* I_G(a) + o(n).
\end{equation}
Combining relations~\eqref{eqNablauntn1} and~\eqref{equntn2''Grad}, we conclude that
$$
\nabla u_n(t_n)=\big(-2c-(h_1+h_2) p_* I_G(a)\big)n + o(n).
$$
On the other hand, $\nabla u_n(t_n)=O(1)$ by Proposition~\ref{lGradientuBounded}. Hence,~\eqref{eqa2} follows.

%\begin{remark}
%Relation~\eqref{equntn2Grad} and Lemma~\ref{lGradientuBounded} actually imply that the remainder in~\eqref{eqNablauntnWithLog} must be $O(1)$. One could get this directly if one expanded $\nabla \Gr_{n-k}(a(n^2-k^2)+\omega_n-\omega_k)$ from~\eqref{eqGraduntn1'} into the Taylor series up to $\nabla\ddot \Gr_{n-k}$. Then, using Proposition~\ref{theoremyn} for $\nabla\dot \Gr_{n-k}(a(n^2-k^2))$ and the fact that $h'(0)=0$, one would see that $\hat\Theta_{1,n}$ is actually $O(1)$. This would yield the remainder in~\eqref{eqNablauntnWithLog} of order $O(1)$ instead of $O(\ln n)$.
%\end{remark} 

%% file: sect7.tex
\section{Proof of the main result: Theorem~$\ref{thMain}$}\label{secProofMainResult}

{\bf 1.} The existence of the limits in~\eqref{eqExistLimitpq} follows from Proposition~\ref{propAsympun}. Let us prove the rest. Denote $b:=\dfrac{h_1+h_2}{h_1}p_*$ and rewrite equalities~\eqref{eqa1} and~\eqref{eqa2} as follows:
\begin{align}
-c+(h_1-2c)a- b h_1 I_F(a) & =0,\label{eqFab} \\
-2c-b h_1 I_G(a) & =0.\label{eqGab}
\end{align}
It is proved in~\cite[Proposition 4.2]{GurTikhRattling1-0} that each of the equations
\begin{align}
-c+(h_1-2c)a_*- h_1 I_F(a_*) & =0,\label{eqFa} \\
-2c-h_1 I_G(a_*) & =0\label{eqGa}
\end{align}
(with respect to $a_*>0$) has the same unique root (and $a_*$ as a function of $h_1/c$ possesses all the properties claimed in Theorem~\ref{thMain}, see also Fig.~\ref{figah1c}). Therefore,  equalities~\eqref{eqFab} and \eqref{eqGab} hold if $a=a_*$ and $b=1$. Let us prove the opposite: if equalities~\eqref{eqFab} and \eqref{eqGab} hold for some $a>0$ and $b\in\bbR$, then $a=a_*$ and $b=1$.

Eliminating $b$ from~\eqref{eqFab} and~\eqref{eqGab}, we find that
\begin{equation}\label{eqMain1}
\dfrac{h_1}{c}=\dfrac{(2a+1)I_G(a)-2I_F(a)}{a I_G(a)}.
\end{equation}
Assume we have proved that the right-hand side in~\eqref{eqMain1} satisfies
\begin{equation}\label{eqMain2}
\dfrac{(2a+1)I_G(a)-2I_F(a)}{a I_G(a)}=-\dfrac{2}{I_G(a)}.
\end{equation}
Then relations~\eqref{eqMain1} and~\eqref{eqMain2} imply that $-2c-h_1 I_G(a)=0$. As we mentioned, by~\cite[Proposition 4.2]{GurTikhRattling1-0}, this yields $a=a_*$. Now, comparing~\eqref{eqGab} with $a=a_*$ and~\eqref{eqGa}, we conclude that $b=1$.

{\bf 2.} Let us prove~\eqref{eqMain2}. We rewrite it as
\begin{equation}\label{eqMain3}
(2a+1)I_G(a)-2I_F(a)+2a=0.
\end{equation}

Making the change of variables $y=\dfrac{1}{\sqrt{a}}\sqrt{\dfrac{1-x}{1+x}}$ in the integral in~\eqref{eqG}, integrating by parts, and using~\eqref{eqfgh1} and~\eqref{eqg0}, we have
\begin{equation}\label{eqMain4}
  I_G(a)=-\int_{0}^{\infty}\left(\dfrac{1-ay^2}{1+ay^2}\right)' g(y)\,dy = \int_{0}^{\infty}\dfrac{1-ay^2}{1+ay^2} h(y)\,dy -\dfrac{1}{2}.
\end{equation}
%Making the same change of variables, integrating by parts and using~\eqref{eqfgh1}, we obtain
%\begin{equation*}
%\begin{aligned}
%  I_F(a)&=-\int_{0}^{\infty}\left(\dfrac{1-ay^2}{1+ay^2}\right)' \dfrac{2ay}{1+ay^2} f(y)\,dy \\
%  & = \int_{0}^{\infty}\dfrac{1-ay^2}{1+ay^2} \left(\dfrac{2ay}{1+ay^2}g(y)+\dfrac{2a(1-ay^2)}{(1+ay^2)^2}f(y)\right)\,dy.
%\end{aligned}
%\end{equation*}
%Eliminating $f(y)$ with the help of~\eqref{eqfgh2} yields
%\begin{equation*}
%\begin{aligned}
%  I_F(a)&= \int_{0}^{\infty}\dfrac{4a(1-ay^2)}{(1+ay^2)^3} (yg(y)+(1-ay^2)h(y))\,dy.
%\end{aligned}
%\end{equation*}
%Noting that $\left(\dfrac{2ay^2}{(1+ay^2)^2}\right)'=\dfrac{4ay(1-ay^2)}{(1+ay^2)^3}$ (and using $g'(y)=h(y)$), we obtain
%\begin{equation}\label{eqMain5}
%\begin{aligned}
%  I_F(a)&= \int_{0}^{\infty}\dfrac{4a(1-ay^2)}{(1+ay^2)^3} \left(yg(y)+(1-ay^2)h(y)\right)\,dy.
%\end{aligned}
%\end{equation}
Making the same change of variables in the integral in~\eqref{eqF} and additionally using~\eqref{eqfgh1}, one can show (see details in~\cite[Section~B]{GurTikhRattling1-0}) that
\begin{equation}\label{eqMain5}
  I_F(a)= \int_{0}^{\infty}\dfrac{2a+1}{1+ay^2} h(y)\,dy -\dfrac{1}{2}.
\end{equation}

Due to~\eqref{eqMain4} and~\eqref{eqMain5}, the left-hand side of relation~\eqref{eqMain3} equals
$$
(2a+1)I_G(a)-2I_F(a)+2a=(2a+1)\left(\dfrac{1}{2}-\int_0^\infty h(y)\,dy\right)=0.
$$

{\bf 3.} It remains to prove the properties of the function $a_*=a_*(h_1/c)$ stated in item 2 of Theorem~\ref{thMain}. In order to prove that the function $a_*$ is infinitely differentiable and decreasing, due to~\eqref{eqGa}, it suffices to show that $I_G(a)$ is infinitely differentiable and $I_G'(a) < 0$ for $a>0$.  The former easily follows from the definition of $I_G$ in~\eqref{eqG} and the properties of the function~$g$ described in Appendix~\ref{secAppendixAsymptoticsFS}. To see the latter, we use~\eqref{eqG} and~\eqref{eqfgh1} to obtain
$$
I_G'(a)=-\dfrac{1}{2a^{3/2}}\int\limits_{-1}^1 h\left(\dfrac{1}{\sqrt{a}}\sqrt{\dfrac{1-x}{1+x}}\right)\sqrt{\dfrac{1-x}{1+x}}\,dx<0.
$$
In order to prove \eqref{AddTh}, it is enough to show that
$$
\mbox{$I_G(a) \to -1$ as $a \to + \infty$} \quad \mbox{and} \quad \mbox{$I_G(a) \to 0$ as $a \to 0$},
$$
which trivially follows from \eqref{eqG}, \eqref{eqg0}, and the properties of the function~$g$ described in Appendix~\ref{secAppendixAsymptoticsFS}.

Theorem~\ref{thMain} is proved.

\begin{figure}[ht]
  \centering
  \includegraphics[width=0.75\linewidth]{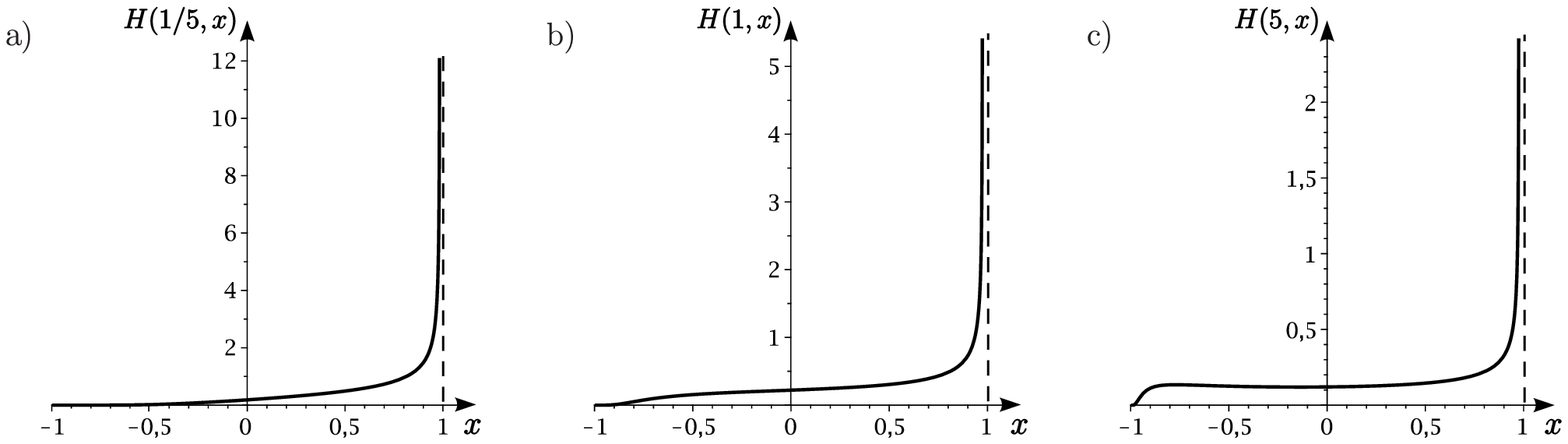}
  \caption{Graphs of function $H(a, \cdot)$ for $a)$ $a = 1/5$, $b)$ $a = 1$, $c)$ $a = 5$.}\label{figH}
\end{figure}

\begin{remark}
Finding $a_*$ from equation~\eqref{eqFa} or~\eqref{eqGa} requires evaluation of the integral $I_F(a)$ or $I_G(a)$ (see~\eqref{eqF} and~\eqref{eqG}) involving the function $f$ or $g$, respectively. In turn, the latter two are defined via the integral of the function $h$, see~\eqref{eqfghft}. From the computational point of view, it would be desirable to avoid evaluating  multiple integrals. In our case, it appears to be possible.
Similarly to~\eqref{eqF} and~\eqref{eqG}, set
$$
H(a,x):=\dfrac{1}{\sqrt{a}{\sqrt{1-x^2}}}\, h\left(\dfrac{1}{\sqrt{a}}\sqrt{\dfrac{1-x}{1+x}}\right),\quad
I_H(a):=\int\limits_{-1}^1 H(a,x)\,dx.
$$
$($see Fig.~$\ref{figH}$$)$, where the function $h$ is defined in~\eqref{eqfghft}. Then~\cite[Proposition 4.2]{GurTikhRattling1-0} implies that the equation
\begin{equation}\label{eqHa}
  (h_1-2c)-h_1 I_H (a_*)=0,\quad a_*>0,
\end{equation}
has a unique root, which coincides with the root of equations~\eqref{eqFa} and~\eqref{eqGa}.
\end{remark}

%% file: appendixA.tex
\appendix

\section{Discrete Green function $\Gr_n(t)$}\label{secAppendixAsymptoticsFS}

In this appendix, we
introduce one of our main tools, namely, the   {\it
discrete Green function}
\begin{equation}\label{eqynIntegralFormula}
\Gr_n(t):=\dfrac{1}{2\pi}\int_{-\pi}^\pi\dfrac{1-e^{-2t(1-\cos\theta)}}{2(1-\cos\theta)}\,e^{i
n\theta}\,d\theta, \quad t\ge 0.
\end{equation}
For convenience, we set $\Gr_n(t)=0$ for $t<0$.
One can check that $\Gr_n\in C^\infty[0,\infty)$ and satisfies
\begin{equation}\label{eqy_nGreen}
\left\{
\begin{aligned}
& \dot \Gr_0=\Delta \Gr_0+1, & & t>0,\\
& \dot \Gr_n=\Delta \Gr_n, & & t>0,\ n\ne0,\\
& \Gr_n(0)=0, & & n\in\mathbb Z.
\end{aligned}
\right.
\end{equation}

  Asymptotic properties of $\Gr_n(t)$ and its derivatives play a central role in our paper. To describe these properties, we will use the functions $h, f, g: \bbR^+ \to \bbR$ given by
\begin{equation}\label{eqfghft}
h(x) := \dfrac{1}{2\sqrt{\pi}}\,e^{-{x^2}/{4}}, \quad
f(x) := 2x\int_x^\infty y^{-2}h(y)\,dy, \quad
g(x) := f'(x).
\end{equation}

\begin{figure}[ht]
  \centering
  \includegraphics[width=0.75\linewidth]{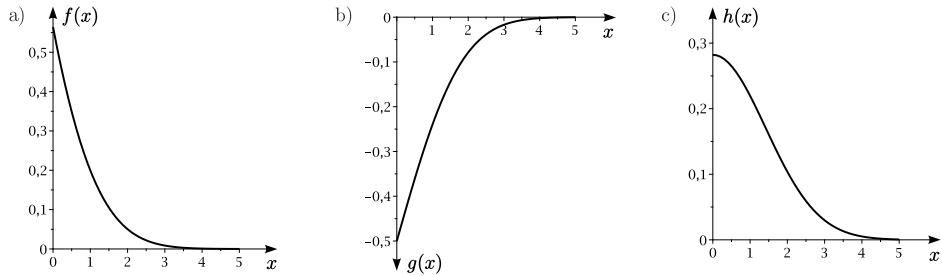}
  \caption{Graphs of functions $f$, $g$, and $h$ }\label{figfgh}
\end{figure}

Note that $f,g,h\in C^{\infty}[0, \infty)$ and decay to zero as $x \to +\infty$, together with all their derivatives, faster than exponentially, see Fig.~\ref{figfgh}. Moreover, they satisfy the relations
\begin{gather}
h(x) = g'(x) = f''(x), \label{eqfgh1}\\
g(0)=-\dfrac{1}{2},\label{eqg0}\\
2h(x) + xg(x) - f(x) = 0. \label{eqfgh2}
\end{gather}

%Consider the functions $r_1, r_2: \bbZ \times \bbR_+ \to \bbR$ satisfying the relations
%$$
%\begin{aligned}
%\Gr_n(t) & = \sqrt{t} f(n/\sqrt{t}) + r_1(n, t),\\
%\nabla \Gr_n(t) & = g(n/\sqrt{t}) +  r_2(n, t).
%\end{aligned}
%$$
%The following result was proved in~\cite{AsympGreenFunc}.
\begin{prp}[see \cite{AsympGreenFunc}]\label{theoremyn}
For any $t_*>0$, there exist constants $C_1,\dots,C_4 > 0$ such that, for all $t \geq
t_*$ and $n = 0, 1, 2, \dots$, we have
$$
\begin{aligned}
 &\left|\Gr_n(t) - \sqrt{t} f\left(\dfrac{n}{\sqrt{t}}\right)\right| \leq \frac{C_1}{\sqrt{t}},  &  &\left|\nabla \Gr_n(t) - g\left(\dfrac{n}{\sqrt{t}}\right)\right|  \leq  \frac{C_2}{\sqrt{t}}, \\
 & \left|\dot{\Gr}_n(t)\right| \leq  \frac{C_3}{\sqrt{t}}, &   &\left|\nabla\dot \Gr_n(t)\right| \leq \frac{C_4}{t}.
\end{aligned}
$$
If $n\le -1$, then $\Gr_n(t)=\Gr_{-n}(t)$ and $\nabla\Gr_n(t)=-\nabla \Gr_{-(n+1)}(t)$.
\end{prp}

{\bf Acknowledgements.}  The work of the first author was
supported by the DFG Heisenberg Programme and the DFG project SFB 910. The work of the second author
was partially supported by Chebyshev Laboratory (Department of Mathematics
and Mechanics, Saint-Petersburg State University) under RF Government grant
11.G34.31.0026, JSC “Gazprom neft”, the Saint-Petersburg State University
research grant 6.38.223.2014, and RFBR 15-01-03797a.